\documentclass[12pt]{article}
\usepackage{amssymb}
\usepackage{amsfonts}
\usepackage{amsmath}

\input{tcilatex}
\begin{document}

\title{Diophantine approximation with improvement of the simultaneous
control of the error and of the denominator }
\author{Abdelmadjid BOUDAOUD \\
Department of Mathematics, \\
Faculty of Mathematics and Computer Sciences, \\
University of M'sila, Algeria\\
Laboratory of Pure and Applied Mathematics (L.M.P.A.)}
\maketitle

\begin{abstract}
In this work we proof the following theorem which is, in addition to some
other lemmas, our main result:

\noindent \textbf{theorem}. Let$\ X=\left\{ \left( x_{1}\text{, }%
t_{1}\right) \text{, }\left( x_{2}\text{, }t_{2}\right) \text{, ..., }\left(
x_{n}\text{, }t_{n}\right) \right\} $ be a finite part of $\mathbb{R}\times 
\mathbb{R}^{\ast +}$, then there exist a finite part $R$ of $\mathbb{R}%
^{\ast +}$ such that for all $\varepsilon >0$ there exists $r\in R$ such
that if $0<\varepsilon \leq r$ then there exist rational numbers $\left( 
\dfrac{p_{i}}{q}\right) _{i=1,2,...,n}$ such that:

\begin{equation}
\left\{ 
\begin{array}{c}
\left\vert x_{i}-\dfrac{p_{i}}{q}\right\vert \leq \varepsilon t_{i} \\ 
\varepsilon q\leq t_{i}%
\end{array}%
\right\vert \text{, }i=1,2,...,n\text{.}  \tag{*}
\end{equation}

\noindent It is clear that the condition $\varepsilon q\leq t_{i}$ for $%
i=1,2,...,n$ is equivalent to $\varepsilon q\leq t=\underset{i=1,2,...,n}{Min%
}$ $\left( t_{i}\right) $.\ Also, we have (*) for all $\varepsilon $
verifying $0<\varepsilon \leq \varepsilon _{0}=\min R$.

The previous theorem is the classical equivalent of the following one which
is formulated in the context of the nonstandard analysis ($\left[ 2\right] $%
, $\left[ 5\right] $, $\left[ 6\right] $, $\left[ 8\right] $).

\noindent \textbf{theorem. }For every positive infinitesimal real $%
\varepsilon $, there exists an unlimited integer $q$\ depending only of $%
\varepsilon $, such that\textit{\ }$\forall ^{st}x\in \mathbb{R}$ $\exists $ 
$p_{x}\in \mathbb{Z}$:

\begin{equation*}
\left\{ 
\begin{array}{ccc}
x & = & \dfrac{p_{x}}{q}+\varepsilon \phi \\ 
\varepsilon q & \cong & 0%
\end{array}%
\text{ .}\right.
\end{equation*}

For this reason,\ to prove the nonstandard version of the main result and to
get its classical version\ we place ourselves in the context of the
nonstandard analysis.
\end{abstract}

\noindent \textbf{1991 Mathematics Subject Classification.} 11J13, 03H05,
26E35.

\noindent \textbf{Key words and phrases.} Diophantine approximation, Farey
series, Nonstandard Analysis.

\section{Introduction, Notations and Rappel}

We dispose in the domain of Diophantine approximation of many results (refer
for example to $\left[ 3\right] $, $\left[ 7\right] $). In the following, we
give as an example, the two most used theorems:

\noindent \textbf{Theorem (Dirichlet) 1.1. }$\left[ \text{7}\right] $\textbf{%
.\ }Suppose that \textit{\ }$x_{1}\mathit{,\ }x_{2}\mathit{,\ }...\mathit{\ ,%
}x_{n}$\textit{\ are }$n$ \textit{real numbers and that }$T>1$\textit{\ is
an integer. Then there exist integers }$q$\textit{,}$p_{1}$\textit{,}$p_{2}$%
\textit{,}$...$\textit{,}$p_{n}$\textit{\ with }

\begin{equation}
\left\{ 
\begin{array}{ccc}
\left\vert x_{i}-\dfrac{p_{i}}{q}\right\vert \leq \dfrac{1}{Tq} &  & \left(
i=1,2,...,n\right) \\ 
1\leq q<T^{n} &  & 
\end{array}%
\text{ .}\right.  \tag{1.1}
\end{equation}

\noindent \textbf{Theorem (Kronecker) 1.2. }$\left[ \text{7}\right] $\textbf{%
.} For any reals $\beta _{1}\mathit{,\ }\beta _{2}\mathit{,\ }...\mathit{\ ,}%
\beta _{n}$ and any $t>0$, the system of inequalities

\begin{equation}
\left\{ 
\begin{array}{c}
\left\vert q\zeta _{1}-p_{1}-\beta _{1}\right\vert <t \\ 
\left\vert q\zeta _{2}-p_{2}-\beta _{2}\right\vert <t \\ 
.................. \\ 
\left\vert q\zeta _{n}-p_{n}-\beta _{n}\right\vert <t%
\end{array}%
\right\vert  \tag{1.2}
\end{equation}%
is solvable in integers $q$, $p_{1}\mathit{,\ }p_{2}\mathit{,\ }...\mathit{\
,}p_{n}$ if and only if $\zeta _{1}\mathit{,\ }\zeta _{2}\mathit{,\ }...%
\mathit{\ ,}\zeta _{n}$ are not rationally dependent. Note that $\zeta _{1}%
\mathit{,\ }\zeta _{2}\mathit{,\ }...\mathit{\ ,}\zeta _{n}$ are said
rationally dependent if there exist integers $r$, $r_{1}\mathit{,\ }r_{2}%
\mathit{,\ }...\mathit{\ ,}r_{n}$ not all zero such that

\begin{equation*}
r_{1}\zeta _{1}+r_{2}\zeta _{2}+...+r_{n}\zeta _{n}=r\text{.}
\end{equation*}%
When we take $\beta _{1}=\beta _{2}=...\mathit{\ }=\beta _{n}=0$, this
theorem is used to approximate the reals $\zeta _{i}$ by using rationals $%
\dfrac{p_{i}}{q}$ to errors smaller than $\dfrac{t}{q}$.

In general, in these results we observe that the simultaneous control
between the error and the common denominator $q$ should be clarified and
specified. This, because the approximation to a given error (which is
generally small) requires a denominator that is generally too big.
Conversely, the approximation with a small denominator might give an error
that is not really small. This question has motivated us to give the
following theorem which is, in addition to some other lemmas, our main
result of this work.\newline

\noindent \textbf{Theorem 1.3.} Let$\ X=\left\{ \left( x_{1}\text{, }%
t_{1}\right) \text{, }\left( x_{2}\text{, }t_{2}\right) \text{, ..., }\left(
x_{n}\text{, }t_{n}\right) \right\} $ be a finite part of $\mathbb{R}\times 
\mathbb{R}^{\ast +}$, then there exist a finite part $R$ of $\mathbb{R}%
^{\ast +}$ such that for all $\varepsilon >0$ there exists $r\in R$ such
that if $0<\varepsilon \leq r$ then there exist rational numbers $\left( 
\dfrac{p_{i}}{q}\right) _{i=1,2,...,n}$ such that:

\begin{equation}
\left\{ 
\begin{array}{c}
\left\vert x_{i}-\dfrac{p_{i}}{q}\right\vert \leq \varepsilon t_{i} \\ 
\varepsilon q\leq t_{i}%
\end{array}%
\right\vert \text{, }i=1,2,...,n\text{.}  \tag{1.3}
\end{equation}%
We note that in (1.3) the condition $\varepsilon q\leq t_{i}$ for $%
i=1,2,...,n$ \ is equivalent to $\varepsilon q\leq t=\underset{i=1,2,...,n}{%
Min}$ $\left( t_{i}\right) $. Also, under the assumption of theorem 1.3, for
all $\varepsilon $ verifying $0<\varepsilon \leq \varepsilon _{0}=\min R$ we
obtain (1.3).

The theorem 1.3 is the classical equivalent of the following theorem
(theorem 1.4.) formulated in the context of the nonstandard analysis.

\noindent \textbf{Theorem 1.4. }For every positive infinitesimal real $%
\varepsilon $, there exists an integer $Q$\ depending only of $\varepsilon $%
, such that\ $\forall ^{st}x\in R$ $\exists $ $P_{x}\in \mathbb{Z}$:

\begin{equation}
\left\{ 
\begin{array}{ccc}
x & = & \dfrac{P_{x}}{Q}+\varepsilon \phi \\ 
\varepsilon Q & \cong & 0%
\end{array}%
\text{ .}\right.  \tag{1.4}
\end{equation}

In the following we make a comparison between our result (theorem 1.3) and
the existing results such as Dirichlet's theorem and Kronecker's theorem.

Our main result is used to approximate at a reduced common denominator $q$
since $\varepsilon q\leq t$ (i.e. $q\leq \dfrac{t}{\varepsilon }$) and at a
different errors since $\left\vert x_{i}-\dfrac{p_{i}}{q}\right\vert \leq
\varepsilon t_{i}$ for $i=1,2,...,n$. In addition, if we take $%
t_{1}=t_{2}=...=t_{n}=t>0$ and $\varepsilon _{0}=Min$ $R$ \ then for every $%
0<\varepsilon $ $\leq \varepsilon _{0}$ there exist integers $q$,$p_{1}$,$%
p_{2}$,$...$,$p_{n}$\ such that

\begin{equation}
\underset{i\in \left\{ 1,2,...,n\right\} }{Max}\left\vert x_{i}-\dfrac{p_{i}%
}{q}\right\vert \leq \varepsilon t\text{ and }q\leq \dfrac{t}{\varepsilon } 
\tag{1.5}
\end{equation}%
i.e., a denominator $q\leq \dfrac{t}{\varepsilon }$ enough for an error not
exceeding $\varepsilon t$.

Look when we use, under the same hypotheses, the Dirichlet's theorem. It may
happen that when we take $\dfrac{1}{T}>\varepsilon t$, the common
denominator $q\geq 1$ is small enough so that the maximum error is strictly
greater than $\varepsilon t$ i.e. $\varepsilon t<\underset{i\in \left\{
1,2,...,n\right\} }{Max}\left\vert x_{i}-\dfrac{p_{i}}{q}\right\vert \leq 
\dfrac{1}{Tq}\leq \dfrac{1}{T}$. In contrast, when we take $T$ satisfying $%
\dfrac{1}{T}\leq \varepsilon t$ \ then we are sure that the maximum error is
smaller than or equal to $\varepsilon t$ i.e. $\underset{i\in \left\{
1,2,...,n\right\} }{Max}\left\vert x_{i}-\dfrac{p_{i}}{q}\right\vert \leq 
\dfrac{1}{Tq}\leq \dfrac{1}{T}\leq \varepsilon t$. But in this case it may
happen that the common denominator $q$, since that $1\leq q<T^{n}$, is very
close to $T^{n}\geq \dfrac{1}{\left( \varepsilon t\right) ^{n}}$ ( $%
q=T^{n}-1\geq $ $\dfrac{1}{\left( \varepsilon t\right) ^{n}}-1$; for
instance). Consequently, to be sure of the realization of the approximation
asked, it is necessary to choose $\dfrac{1}{T}\leq \varepsilon t$ and $q$
can be too big in this case as we have seen.

On his part the Kronecker's theorem is purely existential and don't say
anything on the common denominator.

From the above we can see that the theorem 1.3 ensure the ability to control
the size of $q$ and of the maximum error; especially when $\varepsilon $
(resp. $n$) become small (resp. large). For its proof we place ourselves in
the framework of the nonstandard analysis and we proceed as follows :\newline

\noindent (1) We first show theorem 1.4 (In the sequel noted theorem 2.1.)
by using some lemmas.

\noindent (2) We translate theorem 1.4 by using the Nelson's algorithm.

\subsection{Notations}

\noindent \textbf{i) }For a number $x$ (integer or non) we have the
following usages:

\noindent 1) Abbreviation, $st(x)$ indicates that $x$ is standard; $\forall
^{st}x$ signifies

\noindent $\forall x\left[ st(x)\Longrightarrow ..\right] $.

\noindent 2) $x\cong +\infty $ ( resp. $x\cong 0$) signifies that $x$ is a
positive unlimited (resp. $x$ an infinitesimal). $x\underset{\cong }{>}0$
signifies that $x$ is an infinitesimal real strictly positive.

\noindent 3) $\pounds $(resp. $\phi $) signifies a limited real (resp. an
infinitesimal real) on which one doesn't say anything besides.

\noindent 4) $\left\Vert x\right\Vert $ is the difference, taken positively,
between $x$ and the nearest integer.

\noindent 5) $E\left( x\right) $ (resp. $\left\{ x\right\} $) is the
integral part of $x$ (resp. the fractional part of $x$; that is $\left\{
x\right\} =$ $x-E\left( x\right) $).

\noindent 6) Let $\varepsilon $ be an infinitesimal real, one designates by $%
\varepsilon -galaxie\left( x\right) $ the set $\left\{ y\text{ : }%
y=x+\varepsilon \pounds \right\} $ and by $\varepsilon -halo\left( x\right) $
the set $\left\{ y\text{ : }y=x+\varepsilon \phi \right\} $.

\noindent 7) $x^{0}$ signifies, for $x$ limited, the standard part of $x$.

\noindent \textbf{ii)}

\noindent 8) If $E$ is a given set, $E^{\sigma }$(resp. $\left\vert
E\right\vert $) designates the external set formed, only, by the standard
elements of $E$ (resp. the cardinality of $E$).

\noindent 9) One notes by $\left( x_{1},x_{2},...,x_{n}\right) ^{T}$ the
vector column $\left( 
\begin{array}{c}
x_{1} \\ 
x_{2} \\ 
\vdots \\ 
x_{n}%
\end{array}%
\right) $.

\subsection{\textbf{Rappel}}

\subsubsection{\textbf{Farey series}($\left[ 3\right] $)}

\noindent The Farey series $\mathcal{F}_{N}$ of order $N$ is the ascending
series of irreducible fractions between $0$ and $1$ whose denominators do
not exceed $N$. Thus $\dfrac{h}{k}$ belongs to $\mathcal{F}_{N}$ if

\begin{center}
$0\leq h\leq k\leq N$, $\left( h,k\right) =1$
\end{center}

\noindent the numbers $0$ and $1$ are included in the forms $\dfrac{0}{1}$
and $\dfrac{1}{1}$. If $\dfrac{h}{k}<\dfrac{h^{^{\prime }}}{k^{^{\prime }}}<%
\dfrac{h^{^{\prime \prime }}}{k^{^{\prime \prime }}}$ are three successive
elements of $\mathcal{F}_{N}$ $\left( N>1\right) $, then one has the
following properties:\newline

\noindent 1$^{0}$) $kh^{^{\prime }}-hk^{^{\prime }}=1$.

\noindent 2$^{0}$) $\dfrac{h^{^{\prime }}}{k^{^{\prime }}}=\dfrac{%
h+h^{^{\prime \prime }}}{k+k^{^{\prime \prime }}}$.

\noindent 3$^{0}$) $k+k^{^{\prime }}>N$ and $\dfrac{h}{k}<\dfrac{%
h+h^{^{\prime }}}{k+k^{^{^{\prime }}}}<\dfrac{h^{^{\prime }}}{k^{^{\prime }}}
$.

\noindent 4$^{0}$) If $N>1$, two successive elements of $\mathcal{F}_{N}$
don't have the same denominator.

\noindent 5$^{0}$) Let $\dfrac{h_{1}}{k_{1}}$, $\dfrac{h_{2}}{k_{2}}$ be two
successive elements of $\mathcal{F}_{N}\mathbb{\ }$($N\geq 1$) with $\dfrac{%
h_{1}}{k_{1}}<\dfrac{h_{2}}{k_{2}}$, and let the two following sequences:

\noindent 
\begin{equation}
\left\{ 
\begin{tabular}{l}
$U_{0}=\dfrac{h_{2}}{k_{2}}$, $U_{1}=\dfrac{h_{2}+h_{1}}{k_{2}+k_{1}}$, ... ,%
$U_{i}=\dfrac{h_{2}+ih_{1}}{k_{2}+ik_{1}}$, ... \\ 
$V_{0}=\dfrac{h_{1}}{k_{1}}$, $V_{1}=\dfrac{h_{1}+h_{2}}{k_{1}+k_{2}}$, ..., 
$V_{j}=\dfrac{h_{1}+jh_{2}}{k_{1}+jk_{2}}$, ...%
\end{tabular}%
\text{. }\right.  \tag{1.6}
\end{equation}%
We prove easily that the sequence $\left( U_{i}\right) _{i\in \mathbb{N}}$
(resp. $\left( V_{j}\right) _{j\in \mathbb{N}}$ ) is decreasing (resp.
increasing); besides we have:

\noindent 
\begin{equation}
\left\{ 
\begin{tabular}{l}
$U_{i}-U_{i+1}=\dfrac{1}{\left( k_{2}+ik_{1}\right) \left( k_{2}+\left(
i+1\right) k_{1}\right) }$, $U_{i}-\dfrac{h_{1}}{k_{1}}=\dfrac{1}{%
k_{1}\left( k_{2}+ik_{1}\right) }$ \\ 
$V_{j+1}-V_{j}=\dfrac{1}{\left( k_{1}+jk_{2}\right) \left( k_{1}+\left(
j+1\right) k_{2}\right) }$, $\dfrac{h_{2}}{k_{2}}-V_{j}=\dfrac{1}{%
k_{2}\left( k_{1}+jk_{2}\right) }$%
\end{tabular}%
\text{. }\right.  \tag{1.7}
\end{equation}

\subsubsection{Approximation to the infinitesimal sense of reals}

\noindent \textbf{Theorem 1.5.}\ $\left[ 1\right] $. Let $\xi $ be a real
number. Then for all positive infinitesimal real\textit{\ }$\varepsilon $
there exist a rational number $\dfrac{p}{q}$ and a limited real $l$ such that%
\textit{:}

\begin{equation}
\left\{ 
\begin{array}{ccc}
\xi & = & \dfrac{p_{i}}{q}+\varepsilon l \\ 
\varepsilon q & \cong & 0%
\end{array}%
\right. \text{.}  \tag{1.8}
\end{equation}

\section{Simultaneous approximation to the infinitesimal sense of standard
reals}

We prove in this section the following theorem whose translation by the
algorithm of Nelson gives the theorem 1.3 .\newline

\noindent \textbf{Theorem 2.1. }For every positive infinitesimal real $%
\varepsilon $, there exists \textit{an integer }$Q$\textit{\ depending only
of }$\varepsilon $, such that\textit{\ }$\forall ^{st}x\in \mathbb{R}$ $%
\exists $ $P_{x}\in \mathbb{Z}$:

\begin{equation}
\left\{ 
\begin{array}{ccc}
x & = & \dfrac{P_{x}}{Q}+\varepsilon \phi \\ 
\varepsilon Q & \cong & 0%
\end{array}%
\text{ .}\right.  \tag{2.1}
\end{equation}%
\newline

Let $\varepsilon $ be a positive infinitesimal real. We need to the
following lemmas\newline

\bigskip

\noindent \textbf{Lemma 2.2.} Let $\left( \xi _{1},\xi _{2},...,\xi
_{N}\right) $ a system of real numbers with $N\geq 1$ limited. Then for all
positive infinitesimal real\textit{\ }$\theta $ there are rational numbers $%
\left( \dfrac{p_{i}}{q}\right) _{i=1,2,...,N}$ and limited reals $\left(
l_{i}\right) _{i=1,2,...,N}$ such that for $i=1,2,...,N$\textit{\ :}

\begin{equation}
\left\{ 
\begin{array}{ccc}
\xi _{i} & = & \dfrac{p_{i}}{q}+\theta l_{i} \\ 
\theta q & \cong & 0%
\end{array}%
\right. \text{.}  \tag{2.2}
\end{equation}%
\textbf{Proof.}\ Consider, for every $n\in \mathbb{N}^{\ast }$, the formula:

\begin{center}
$%
\begin{array}{cc}
B\left( n\right) = & 
\begin{array}{c}
\text{"}\forall \text{ }\left( \xi _{1},\xi _{2},...,\xi _{n}\right) \in 
\mathbb{R}^{n}\text{ with }n\geq 1\text{ and }\forall \theta \underset{\cong 
}{>}0\text{ }\exists \left( \dfrac{P_{i}}{Q}\right) _{i=1,2,...,n} \\ 
\text{such that for every }i\in \left\{ 1,2,...,n\right\} \text{ : }\left\{ 
\begin{array}{c}
x_{i}-\dfrac{P_{i}}{Q}=\theta \pounds  \\ 
\theta Q\cong 0%
\end{array}%
\right. \text{ \ "}%
\end{array}%
\text{.}%
\end{array}%
$
\end{center}

By theorem 1.5, we have $B\left( 1\right) $. Suppose, for $1\leq n$ a
standard integer, $B\left( n\right) $ and prove $B\left( n+1\right) $. Let $%
\left( \xi _{1},\xi _{2},...,\xi _{n},\xi _{n+1}\right) \in \mathbb{R}^{n+1}$
and let $\theta \underset{\cong }{>}0$, then by $B\left( n\right) $ there
are rational numbers $\left( \dfrac{p_{i}}{q}\right) _{i=1,2,...,n}$ such
that

\begin{equation}
\left\{ 
\begin{array}{ccc}
\xi _{1} & = & \dfrac{p_{1}}{q}+\theta \pounds  \\ 
\xi _{2} & = & \dfrac{p_{2}}{q}+\theta \pounds  \\ 
\vdots & = & \vdots \\ 
\xi _{n} & = & \dfrac{p_{n}}{q}+\theta \pounds 
\end{array}%
\right.  \tag{2.3}
\end{equation}%
where $\theta q\cong 0$. Now, since $\theta q\cong 0$, the application of
theorem 1.5\textbf{\ } implies $q\xi _{n+1}=\dfrac{p_{n+1}}{q_{n+1}}+\left(
\theta q\right) \pounds $, $\left( \theta q\right) q_{n+1}\cong 0$. Hence

\begin{equation}
\xi _{n+1}=\dfrac{p_{n+1}}{qq_{n+1}}+\theta \pounds \text{, }\theta
qq_{n+1}\cong 0\text{.}  \tag{2.4}
\end{equation}%
We deduct from (2.3) and (2.4) that:

\begin{equation*}
\left\{ 
\begin{array}{ccccc}
\xi _{1} & = & \dfrac{p_{1}q_{n+1}}{qq_{n+1}}+\theta \pounds  & = & \dfrac{%
P_{1}}{Q}+\theta \pounds  \\ 
\xi _{2} & = & \dfrac{p_{2}q_{n+1}}{qq_{n+1}}+\theta \pounds  & = & \dfrac{%
P_{2}}{Q}+\theta \pounds  \\ 
\vdots & = & \vdots & = & \vdots \\ 
\xi _{n} & = & \dfrac{p_{n}q_{n+1}}{qq_{n+1}}+\theta \pounds  & = & \dfrac{%
P_{n}}{Q}+\theta \pounds  \\ 
\xi _{n+1} & = & \dfrac{p_{n+1}}{qq_{n+1}}+\theta \pounds  & = & \dfrac{%
P_{n+1}}{Q}+\theta \pounds 
\end{array}%
\right.
\end{equation*}%
where, from (2.4), $\theta Q=\theta qq_{n+1}\cong 0$. Consequently $B\left(
n+1\right) $. Therefore, by the external recurrence principle, we have $%
\forall ^{st}n\geq 1$ $B\left( n\right) $.\newline

\noindent \textbf{Lemma 2.3.}\textit{\ }Let $E$\ be a given set. \ For all
integer $\omega $\ $\cong +\infty $, there is a finite subset $F\subset E$
containing all standard elements of \ $E$\ (i.e. $E^{\sigma }\subset F$\ )
and whose cardinal is strictly inferior to $\omega $\ ($\left\vert
F\right\vert <\omega $).\newline

\noindent \textbf{Proof.} Let $\omega $\textit{\ }$\cong +\infty $. Let $%
B\left( F,z\right) $ be the internal formula:\ "$F\subset E$, $\left\vert
F\right\vert <\omega $, $z\in F$ ". Let $Z\subset E$ be a standard \ finite
part. Then there exists a finite part $F\subset E$ with $\left\vert
F\right\vert <\omega $ such that every element $z$ of $Z$ belongs to $F$,
i.e. we have $B\left( F,z\right) $. Indeed it suffices to take $F=Z$.
Therefore, the principle of idealization (I) asserts the existence of a
finite part $F\subset E$ with $\left\vert F\right\vert <\omega $ such that
any standard element of $L$ belongs to $F$.\newline

\noindent \textbf{Lemma 2.4. }Let $\lambda \cong +\infty $\ be a real number
such that $\sqrt{\varepsilon }\lambda \cong 0$. Let $F_{M}$\ be the Farey
sequence of order $M=E\left( \dfrac{\lambda }{\sqrt{\varepsilon }}\right) $.
If $\dfrac{p_{1}}{q_{1}}$,\ $\dfrac{p_{2}}{q_{2}}$\ are two elements of $%
F_{M}$\ such that $q_{1}\simeq +\infty $, $q_{2}\simeq +\infty $\ and $\left[
\dfrac{p_{1}}{q_{1}}\text{ , }\dfrac{p_{2}}{q_{2}}\right] $\ doesn't contain
any standard rational number (in this case $\dfrac{p_{1}}{q_{1}}\cong \dfrac{%
p_{2}}{q_{2}}$). Then there exist a finite sequence of irreducible rational
numbers $\left( \dfrac{l_{i}}{m_{i}}\right) _{i=1,2,...,g}$such that:

\begin{equation*}
\dfrac{p_{1}}{q_{1}}=\dfrac{l_{1}}{m_{1}}<\dfrac{l_{2}}{m_{2}}<...<\dfrac{%
l_{g}}{m_{g}}=\dfrac{p_{2}}{q_{2}}
\end{equation*}%
where $\dfrac{l_{i+1}}{m_{i+1}}-\dfrac{l_{i}}{m_{i}}=\varepsilon \phi $\ for 
$i=1,2,...,g-1$. Besides for $i=1,2,...,g$\ we have $\varepsilon m_{i}\cong
0 $\ and $m_{i}\cong +\infty $.\newline

\noindent Proof. Let us consider the case where $\dfrac{p_{2}}{q_{2}}-\dfrac{%
p_{1}}{q_{1}}$ is not of $\varepsilon \phi $ form; otherwise the lemma is
proved. Let $\left( \dfrac{t_{i}}{\gamma _{i}}\right) _{i=1,2,...,r}$ be the
elements of $\mathcal{F}_{M}$ such that

\begin{equation*}
\dfrac{p_{1}}{q_{1}}=\dfrac{t_{1}}{\gamma _{1}}<\dfrac{t_{2}}{\gamma _{2}}%
<...<\dfrac{t_{r}}{\gamma _{r}}=\dfrac{p_{2}}{q_{2}}\text{.}
\end{equation*}%
Let $i_{0}\in \left\{ 1,2,...,r-1\right\} $ such that $\dfrac{t_{i_{0}+1}}{%
\gamma _{i_{0}+1}}-\dfrac{t_{i_{0}}}{\gamma _{i_{0}}}$ is not of $%
\varepsilon \phi $ form, because if a such $i_{0}$ does not exist the lemma
is proved. From the properties of $\mathcal{F}_{M}$ ( \textbf{1.2.1}), $%
\gamma _{i_{0}+1}$ and $\gamma _{i_{0}}$ cannot be equal. Then there are two
cases:

\noindent \textbf{A)} $\gamma _{i_{0}+1}>\gamma _{i_{0}}$ : Let us take, in
this case, $g_{0}\cong +\infty $ an integer such that $\dfrac{g_{0}}{\gamma
_{i_{0}}}\cong 0$ ( the existence of $g_{0}$ is assured by Robinson's
lemma). Let $X=E\left( \dfrac{g_{0}}{\varepsilon \gamma _{i_{0}}}\right) $
and

\begin{equation*}
H=\left\{ \dfrac{t_{i_{0}}}{\gamma _{i_{0}}},U_{p},U_{p-1},...,U_{0}\right\}
\end{equation*}%
where $p=E\left( \dfrac{X-\gamma _{i_{0}+1}}{\gamma _{i_{0}}}\right) $ and $%
U_{i}=\dfrac{t_{i_{0}+1}+i.t_{i_{0}}}{\gamma _{i_{0}+1}+i.\gamma _{i_{0}}}$ (%
$i=0,1,...,p-1,p$). Now we prove that : $p$ is an unlimited integer, the
product of the denominator of every element of $H$ by $\varepsilon $ is an
infinitesimal and the distance between two successive elements of $H$ is of
the $\varepsilon \phi $ form\textit{.}

Indeed, we have $X=E\left( \dfrac{g_{0}}{\varepsilon \gamma _{i_{0}}}\right)
=\dfrac{g_{0}}{\varepsilon \gamma _{i_{0}}}-\rho _{X}$ where $\rho _{X}\in %
\left[ 0,1\right[ $.

\noindent $%
\begin{array}{lll}
\dfrac{X-\gamma _{i_{0}+1}}{\gamma _{i_{0}}} & = & \dfrac{g_{0}}{\varepsilon
\gamma _{i_{0}}\gamma _{i_{0}}}-\dfrac{\rho _{X}}{\gamma _{i_{0}}}-\dfrac{%
\gamma _{i_{0}+1}}{\gamma _{i_{0}}}\text{ } \\ 
&  &  \\ 
& = & \dfrac{g_{0}-\varepsilon \gamma _{i_{0}}\rho _{X}-\varepsilon \gamma
_{i_{0}}\gamma _{i_{0}+1}}{\varepsilon \gamma _{i_{0}}\gamma _{i_{0}}}\text{.%
}%
\end{array}%
$

\noindent Since $\varepsilon \gamma _{i_{0}}\rho _{X}\cong 0$, $\varepsilon
\gamma _{i_{0}}\gamma _{i_{0}+1}$ is a limited real number otherwise $\dfrac{%
t_{i_{0}+1}}{\gamma _{i_{0}+1}}-\dfrac{t_{i_{0}}}{\gamma _{i_{0}}}=\dfrac{1}{%
\gamma _{i_{0}}\gamma _{i_{0}+1}}=\varepsilon \phi $ what contradicts the
supposition. Then $g_{0}-\varepsilon \gamma _{i_{0}}\rho _{X}-\varepsilon
\gamma _{i_{0}}\gamma _{i_{0}+1}$ is a positive unlimited real. On the other
hand $\varepsilon \gamma _{i_{0}}\gamma _{i_{0}}$ is limited; then $\dfrac{%
X-\gamma _{i_{0}+1}}{\gamma _{i_{0}}}$ is a positive unlimited real,
therefore $p$ is also. The greatest denominator in $H$ is $\gamma
_{i_{0}+1}+p\gamma _{i_{0}}$ where $p=\dfrac{g_{0}}{\varepsilon \gamma
_{i_{0}}\gamma _{i_{0}}}-\dfrac{\rho _{X}}{\gamma _{i_{0}}}-\dfrac{\gamma
_{i_{0}+1}}{\gamma _{i_{0}}}-\rho $ with $\rho \in \left[ 0,1\right[ $.

\noindent $%
\begin{array}{lll}
\varepsilon \left( \gamma _{i_{0}+1}+p\gamma _{i_{0}}\right) & = & 
\varepsilon \left( \gamma _{i_{0}+1}+\left( \dfrac{g_{0}}{\varepsilon \gamma
_{i_{0}}\gamma _{i_{0}}}-\dfrac{\rho _{X}}{\gamma _{i_{0}}}-\dfrac{\gamma
_{i_{0}+1}}{\gamma _{i_{0}}}-\rho \right) \gamma _{i_{0}}\right) \\ 
&  &  \\ 
& = & \varepsilon \gamma _{i_{0}+1}+\dfrac{g_{0}}{\gamma _{i_{0}}}%
-\varepsilon \rho _{X}-\varepsilon \gamma _{i_{0}+1}-\varepsilon \rho \gamma
_{i_{0}}\cong 0\text{ .}%
\end{array}%
$

\noindent Hence the product of the denominator of every element of $H$ by $%
\varepsilon $ is an infinitesimal. It remains to prove that the distance
between two elements of $H$ is of the $\varepsilon \phi $ form; Indeed: Let $%
i\in \left\{ 0,1,...,p-1\right\} $, from ($1.7$) we have

\begin{center}
$U_{i}-U_{i+1}=\dfrac{1}{\left( \gamma _{i_{0}+1}+i.\gamma _{i_{0}}\right)
\left( \gamma _{i_{0}+1}+\left( i+1\right) .\gamma _{i_{0}}\right) }$ .
\end{center}

\noindent By hypothesis we have $\gamma _{i_{0}+1}>\gamma _{i_{0}}$, then of
properties of Farey's series (\textbf{1.2.1)}) $2\gamma _{i_{0}+1}>\gamma
_{i_{0}+1}+\gamma _{i_{0}}>M$, then $\gamma _{i_{0}+1}>\dfrac{M}{2}$.

\noindent Let $d_{i}=\varepsilon \left( \gamma _{i_{0}+1}+i.\gamma
_{i_{0}}\right) \left( \gamma _{i_{0}+1}+\left( i+1\right) .\gamma
_{i_{0}}\right) $.

\noindent Seen \ that$\ \ \left( \gamma _{i_{0}+1}\right) ^{2}>\left( \dfrac{%
M}{2}\right) ^{2}$, $d_{i}$ is unlimited, therefore $U_{i}-U_{i+1}=%
\varepsilon \phi $. To finish the proof, we have of ($1.7$):

\begin{equation*}
U_{p}-\dfrac{t_{i_{0}}}{\gamma _{i_{0}}}=\dfrac{1}{\left( \gamma
_{i_{0}+1}+p\gamma _{i_{0}}\right) \gamma _{i_{0}}}\text{.}
\end{equation*}

\noindent Let $d_{p}=\varepsilon \left( \gamma _{i_{0}+1}+p.\gamma
_{i_{0}}\right) \gamma _{i_{0}}$, after the replacement by the value of $p$,
we obtain

\begin{equation*}
d_{p}=\varepsilon \gamma _{i_{0}+1}\gamma _{i_{0}}+g_{0}-\varepsilon \rho
_{X}\gamma _{i_{0}}-\varepsilon \gamma _{i_{0}+1}\gamma _{i_{0}}-\varepsilon
\rho \gamma _{i_{0}}\gamma _{i_{0}}\text{.}
\end{equation*}

\noindent Since $\varepsilon \gamma _{i_{0}}\cong 0$, $\varepsilon \gamma
_{i_{0}}\gamma _{i_{0}}$ is limited, then $d_{p}$ is unlimited; hence

\begin{equation*}
U_{p}-\dfrac{t_{i_{0}}}{\gamma _{i_{0}}}=\varepsilon \phi \text{.}
\end{equation*}%
Thus, we end what we perceived.\newline

\noindent \textbf{B)} $\gamma _{i_{0}}>\gamma _{i_{0}+1}$: Let us take, in
this case, $g_{1}\cong +\infty $ an integer such that $\dfrac{g_{1}}{\gamma
_{i_{0}+1}}\cong 0$ (the existence of $g_{1}$ is assured by Robinson's
lemma). Let $\widetilde{X}=E\left( \dfrac{g_{1}}{\varepsilon \gamma
_{i_{0}+1}}\right) $ and

\begin{equation*}
\widetilde{H}=\left\{ V_{0},V_{1},...,V_{p^{^{\prime }}-1},V_{p^{^{\prime
}}},\dfrac{t_{i_{0}+1}}{\gamma _{i_{0}+1}}\right\}
\end{equation*}%
where $p^{^{\prime }}=E\left( \dfrac{\widetilde{X}-\gamma _{i_{0}}}{\gamma
_{i_{0}+1}}\right) $ and $V_{j}=\dfrac{t_{i_{0}}+j.t_{i_{0}+1}}{\gamma
_{i_{0}}+j.\gamma _{i_{0}+1}}$ ( $j=0,1,...,p^{^{\prime }}-1,p^{^{\prime }}$%
). Since the symmetry of this case with the case \textbf{A}) we prove, as in
the case of $H$, that $p^{^{\prime }}$ is an unlimited integer, the product
of the denominator of every element of $\widetilde{H}$ by $\varepsilon $ is
an infinitesimal and the distance between two successive elements of $%
\widetilde{H}$ is of the $\varepsilon \phi $ form\textit{.}

Thus the elements of $H$ (or of $\widetilde{H}$ ) form a subdivision of the
interval $\left[ \dfrac{t_{i_{0}}}{\gamma _{i_{0}}},\dfrac{t_{i_{0}+1}}{%
\gamma _{i_{0}+1}}\right] $. For the other intervals $\left[ \dfrac{t_{i}}{%
\gamma _{i}},\dfrac{t_{i+1}}{\gamma _{i+1}}\right] _{i\in \left\{
1,2,...,r-1\right\} -\left\{ i_{0}\right\} }$ which don't have a length of $%
\varepsilon \phi $ form we do the same construction as we did with$\left[ 
\dfrac{t_{i_{0}}}{\gamma _{i_{0}}},\dfrac{t_{i_{0}+1}}{\gamma _{i_{0}+1}}%
\right] $.

By regrouping rational numbers which subdivide intervals $\left[ \dfrac{t_{i}%
}{\gamma _{i}},\dfrac{t_{i+1}}{\gamma _{i+1}}\right] $ ($i\in \left\{
1,2,...,r-1\right\} $) not having a length of the $\varepsilon \phi $ form
and the rationals which are borders of intervals having a length of the $%
\varepsilon \phi $ form, we obtain the finite sequence $\left( \dfrac{l_{i}}{%
m_{i}}\right) _{i=1,2,...,g}$ . The irreducibility of the elements of the
sequence $\left( \dfrac{l_{i}}{m_{i}}\right) _{i=1,2,...,g}$ results from
properties of Farey's series.$\square $\newline

\noindent \textbf{Lemma 2.5.} Let $\xi \in \left[ 0\text{ , }1\right] $ be a
real, if $\xi $ is not in the $\varepsilon $-galaxie of a standard rational
number then there exists two irreducible rational numbers $\dfrac{h_{1}}{%
k_{1}},\dfrac{h_{2}}{k_{2}}$ of the interval $\left[ 0,1\right] $ such that

\begin{equation*}
\xi \in \left[ \dfrac{h_{1}}{k_{1}},\dfrac{h_{2}}{k_{2}}\right] ,\text{ }%
k_{1}\cong +\infty ,\text{ }k_{2}\cong +\infty ,\text{ }\varepsilon
k_{1}\cong \varepsilon k_{2}\cong 0\text{ and }\dfrac{h_{2}}{k_{2}}-\dfrac{%
h_{1}}{k_{1}}=\varepsilon \phi \text{.}
\end{equation*}

\noindent \textbf{Proof.} Let us take, as in the lemma 2.4, a positive
unlimited real number $\lambda $ such that $\sqrt{\varepsilon }\lambda \cong
0$ and let $F_{M}$\textit{\ be }the Farey sequence of order\textit{\ }$%
M=E\left( \dfrac{\lambda }{\sqrt{\varepsilon }}\right) $. Let $\dfrac{p_{1}}{%
q_{1}}$, $\dfrac{p_{2}}{q_{2}}$ be two successive elements of $F_{M}$ such
that $\xi \in \left[ \dfrac{p_{1}}{q_{1}},\dfrac{p_{2}}{q_{2}}\right] $. Two
cases are distinguished:

\noindent \textbf{A) }Nor $\dfrac{p_{1}}{q_{1}}$ nor $\dfrac{p_{2}}{q_{2}}$
is a standard rational : In this case by applying the lemma 2.4, we obtain
two irreducible rationals $\dfrac{l_{i_{0}}}{m_{i_{0}}}$ and $\dfrac{%
l_{i_{0}+1}}{m_{i_{0}+1}}$ such that $\xi \in \left[ \dfrac{l_{i_{0}}}{%
m_{i_{0}}}\text{,}\dfrac{l_{i_{0}+1}}{m_{i_{0}+1}}\right] $, $m_{i_{0}}\cong
+\infty $, $m_{i_{0}+1}\cong +\infty $, $\varepsilon m_{i_{0}}\cong
\varepsilon m_{i_{0}+1}\cong 0$, $\dfrac{l_{i_{0}+1}}{m_{i_{0}+1}}-\dfrac{%
l_{i_{0}}}{m_{i_{0}}}=\varepsilon \phi $. Hence the lemma is proved by
taking $\dfrac{l_{i_{0}}}{m_{i_{0}}}$ for $\dfrac{h_{1}}{k_{1}}$ and $\dfrac{%
l_{i_{0}+1}}{m_{i_{0}+1}}$ for $\dfrac{h_{2}}{k_{2}}$.\newline

\noindent \textbf{B) }$\dfrac{p_{1}}{q_{1}}$ or $\dfrac{p_{2}}{q_{2}}$ is
standard (cannot be both at the same time standard). Let us suppose that $%
\dfrac{p_{1}}{q_{1}}$ is standard (the other case, seen the symmetry, can be
treated by the same way.). Then $\xi -\dfrac{p_{1}}{q_{1}}=\varepsilon w$
where $w\cong +\infty $. Let us put $L=E\left( 2/\left( \xi -\dfrac{p_{1}}{%
q_{1}}\right) \right) $ then $\varepsilon L\cong 0$ and $\dfrac{p_{1}}{q_{1}}%
+\dfrac{1}{L}<\xi $. Let $\dfrac{l}{m}$ be the reduced form of $\dfrac{p_{1}%
}{q_{1}}+\dfrac{1}{L}$, then $\varepsilon m\cong 0$ because $m\leq Lq_{1}$
and $q_{1}$ is a standard. $m>M$ because $\dfrac{l}{m}$ is not an element of 
$\mathcal{F}_{M}$. Therefore $\varepsilon m^{2}$ is an unlimited because $%
\varepsilon m^{2}>\varepsilon M^{2}$ and $\varepsilon M^{2}$ is an
unlimited. This means that $m$ is of the $E\left( \dfrac{\lambda ^{^{\prime
}}}{\sqrt{\varepsilon }}\right) $ form where $\lambda ^{^{\prime }}$ is a
positive unlimited real verifying $\sqrt{\varepsilon }\lambda ^{^{\prime
}}\cong 0$. Now if we consider $\mathcal{F}_{m}$, then $\xi \in \left[ 
\dfrac{p_{1}^{^{\prime }}}{q_{1}^{^{\prime }}},\dfrac{p_{2}^{^{\prime }}}{%
q_{2}^{^{\prime }}}\right] $ where $\dfrac{p_{1}^{^{\prime }}}{%
q_{1}^{^{\prime }}}$ and $\dfrac{p_{2}^{^{\prime }}}{q_{2}^{^{\prime }}}$
are two successive non standard elements of $\mathcal{F}_{m}$. Thus the case 
\textbf{B)} comes back itself to the case \textbf{A)}, therefore the
proposition is also proved for this case.\textbf{\ }$\square $\newline

\noindent \textbf{Remark.} We easily see that this proof is also a proof for
the theorem 1.5.\newline

Let $\gamma $ be a positive unlimited real such that $\varepsilon .\gamma
\simeq 0$, then\newline

\noindent \textbf{Lemma 2.6. }There exists a finite set

\begin{equation}
S=\left\{ l_{1},l_{2},...,l_{n}\right\} \subset \left[ 0,1\right]  \tag{2.5}
\end{equation}%
containing all standard elements of $\left[ 0,1\right] $ such that $%
\left\vert l_{i+1}-l_{i}\right\vert \geq \varepsilon \gamma $ for $i\in
\left\{ 1,2,...,n-1\right\} $.\newline

\noindent \textbf{Proof.} Let $B\left( S,z\right) $ be the internal
formula:\ "$S\subset \left[ 0,1\right] $ is finite, $z\in S$ \& $\forall
\left( x_{1}\text{ , }x_{2}\right) \in S\times S$ $\left( \left\vert
x_{1}-x_{2}\right\vert \geq \varepsilon \gamma \right) $". Let $Z\subset %
\left[ 0,1\right] $ be a standard \ finite part. Then there exists a finite
part $S\subset \left[ 0,1\right] $ such that every element $z$ of $Z$
belongs to $S$ and $\forall \left( x_{1}\text{ , }x_{2}\right) \in S\times S$
$\left( \left\vert x_{1}-x_{2}\right\vert \geq \varepsilon \gamma \right) $,
i.e. we have $B\left( S,z\right) $. Indeed it suffices to take $S=Z$.
Therefore, the principle of idealization (I) asserts the existence of a
finite part $S\subset \left[ 0,1\right] $ such that any standard element of $%
\left[ 0,1\right] $ belongs to $S$ and $\forall \left( x_{1}\text{ , }%
x_{2}\right) \in S\times S$ $\left( \left\vert x_{1}-x_{2}\right\vert \geq
\varepsilon \gamma \right) $. Put $S=\left\{ l_{1},l_{2},...,l_{n}\right\} $%
, where $\left\vert l_{i+1}-l_{i}\right\vert \geq \varepsilon \gamma $ for $%
i\in \left\{ 1,2,...,n-1\right\} $ and any standard element of $\left[ 0,1%
\right] $ belongs to $S$.\newline

\noindent \textbf{Corollary 2.7.} For every element $l_{i}$ of $S$ ($S$%
\textit{\ }is the set that has been constructed in the lemma 2.6 ) we have\
only one of the two cases:

\noindent 1)\textit{\ }$l_{i}$\textit{\ }is a standard rational number.

\noindent 2) $l_{i}$ is outside of $\varepsilon -$galaxies of all standard
rational number.\newline

\noindent \textbf{Proof.} Let $l_{i}\in S$, then

\noindent 1) $l_{i}$ can be a standard rational because $S$ contains all
standard elements of $[0,1]$.

\noindent 2) $l_{i}$ is not a standard rational then $l_{i}$ is not in the $%
\varepsilon -$galaxy of any standard rational. Indeed, suppose that $l_{i}=%
\dfrac{p}{q}+\varepsilon \pounds $ ($\pounds \neq 0$), where $\dfrac{p}{q}$
is standard. Then $l_{i}$ and $\dfrac{p}{q}$ are elements of $S$ with $%
\left\vert l_{i}-\dfrac{p}{q}\right\vert =\left\vert \varepsilon \pounds %
\right\vert <\varepsilon \gamma $ which contradicts lemma 2.6 .\newline

\noindent \textbf{Lemma 2.8. }For every standard integer $n\geq 1$. The real
numbers $x_{i}$\ of all system $\left\{ x_{1},x_{2},...,x_{n}\right\}
\subset S$\ ($S$\ is the set that has been constructed in the lemma 2.6.)
are approximated by rational numbers $\left( \dfrac{P_{i}}{Q}\right)
_{i=1,2,...,n}$\ to $\varepsilon \phi $\ near with $\varepsilon Q$\ $\cong 0$%
. that is to say:

\begin{center}
\begin{equation}
\left\{ 
\begin{array}{ccc}
x_{i} & = & \dfrac{P_{i}}{Q}+\varepsilon \phi \\ 
\varepsilon Q & \cong & 0%
\end{array}%
\right. ;i=1,2,...,n\text{.}  \tag{2.6}
\end{equation}
\end{center}

\noindent \textbf{Proof. }Consider the formula:

\begin{center}
$A(n)\equiv $ $"$ $\forall \left\{ x_{1},x_{2},...,x_{n}\right\} \subset S$ $%
\exists $ $\left( \dfrac{P_{i}}{Q}\right) _{i=1,2,...,n}$ such that:$\left\{ 
\begin{array}{ccc}
x_{i} & = & \dfrac{P_{i}}{Q}+\varepsilon \phi \\ 
\varepsilon Q & \cong & 0%
\end{array}%
\right. $ ; $i=1,2,...,n$ \ $"$.
\end{center}

\noindent According to the corollary 2.7, a real $x$ of $S$ is a standard
rational or is outside of $\varepsilon -$galaxies of standard rationals. In
addition, according to lemma 2.5, if $x$ is not in the $\varepsilon -$galaxy
of a rational standard, $x$ is written in the form $\left\{ 
\begin{array}{ccc}
x & = & \dfrac{P}{Q}+\varepsilon \phi \\ 
\varepsilon Q & \cong & 0%
\end{array}%
\right. $. Then in all cases $x$ is written in the form $\left\{ 
\begin{array}{ccc}
x & = & \dfrac{P}{Q}+\varepsilon \phi \\ 
\varepsilon Q & \cong & 0%
\end{array}%
\right. $. Consequently we have $A\left( 1\right) $.

\noindent Suppose $A\left( n\right) $, for a standard integer $n$, and prove 
$A\left( n+1\right) $.

\noindent Let $\left( x_{1},x_{2},...,x_{n},x_{n+1}\right) \subset S$. Since 
$A$ is verified for $n$ we have

\begin{center}
\begin{equation}
\left\{ 
\begin{array}{ccc}
x_{i}=\dfrac{p_{i}}{q}+\varepsilon \phi & ; & i=1,2,...,n \\ 
\varepsilon q\cong 0 &  & 
\end{array}%
\text{.}\right.  \tag{2.7}
\end{equation}
\end{center}

\noindent If $x_{n+1}=\dfrac{h_{1}}{k_{1}}$ is standard, then because $k_{1}$
is standard and of $(2.7)$ we have

\begin{equation}
\left\{ 
\begin{array}{ccc}
x_{i}=\dfrac{p_{i}k_{1}}{qk_{1}}+\varepsilon \phi =\dfrac{P_{i}}{Q}%
+\varepsilon \phi & ; & i=1,2,...,n \\ 
x_{n+1}=\dfrac{h_{1}q}{k_{1}q}+\varepsilon .0=\dfrac{P_{n+1}}{Q}+\varepsilon
\phi &  &  \\ 
\varepsilon Q=\varepsilon qk_{1}\cong 0 &  & 
\end{array}%
\text{.}\right.  \tag{2.8}
\end{equation}%
Let us look at the case where $x_{n+1}$ is not a rational standard. In this
case the application of the theorem 1.5 to the real $qx_{n+1}$ with the
infinitesimal $\varepsilon q$ implies:

\begin{center}
$\left\{ 
\begin{array}{c}
\begin{array}{ccc}
qx_{n+1} & = & \dfrac{M}{N}+\left( \varepsilon q\right) a%
\end{array}
\\ 
\begin{array}{ccc}
\left( \varepsilon q\right) N & \cong & 0%
\end{array}%
\end{array}%
\right. $
\end{center}

\noindent where $a$ is limited. If $a\cong 0$, then from this and $\left(
2.7\right) $ :

\begin{equation}
\left\{ 
\begin{array}{ccc}
x_{i}=\dfrac{p_{i}N}{qN}+\varepsilon \phi =\dfrac{P_{i}}{Q}+\varepsilon \phi
& ; & i=1,2,...,n \\ 
x_{n+1}=\dfrac{M}{qN}+\varepsilon a=\dfrac{P_{n+1}}{Q}+\varepsilon \phi &  & 
\\ 
\varepsilon Q=\varepsilon qN\cong 0 &  & 
\end{array}%
\text{.}\right.  \tag{2.9}
\end{equation}%
Let us look at the case where $a$ is appreciable. Suppose $a>0$, then

\begin{equation}
\left\{ 
\begin{array}{ccc}
x_{i}=\dfrac{Np_{i}}{Nq}+\varepsilon \phi & ; & i=1,2,...,n \\ 
x_{n+1}=\dfrac{M}{Nq}+\varepsilon a &  &  \\ 
\varepsilon Nq\cong 0 &  & 
\end{array}%
\text{.}\right.  \tag{2.10}
\end{equation}

\noindent The reduced form of $\dfrac{M}{Nq}$ cannot be a rational standard.
Otherwise, $x_{n+1}$ and $\dfrac{M}{Nq}$ become two elements of $S$ such
that the separating distance between them, is of the $\varepsilon a$ form.
What, according to lemma 2.6, is not true for two elements of $S$; for the
same reason $x_{n+1}$ cannot be in the $\varepsilon -$galaxy of a standard
rational. According to the lemma 2.5:%
\begin{equation}
\left\{ 
\begin{array}{lll}
x_{n+1}=\dfrac{h_{1}}{k_{1}}+\varepsilon \phi _{1} & = & \dfrac{h_{2}}{k_{2}}%
-\varepsilon \phi _{2} \\ 
\varepsilon k_{1}\cong \varepsilon k_{2}\cong 0 & ; & k_{1}\cong k_{2}\cong
+\infty%
\end{array}%
\right.  \tag{2.11}
\end{equation}%
Where $\phi _{1}\geq 0$ and $\phi _{2}\geq 0$ are two infinitesimal reals
and $\dfrac{h_{1}}{k_{1}}$, $\dfrac{h_{2}}{k_{2}}$ are irreducibles. Let $%
\xi $ the element of $S$ succeeding immediately $x_{n+1}$ in $S$ ($x_{n+1}<$ 
$\xi $). Then by lemma 2.6 :

\begin{equation*}
\xi -x_{n+1}=\varepsilon \omega \cong 0\text{, }\omega \geq \gamma \text{.}
\end{equation*}%
The real number $\dfrac{x_{n+1}+\xi }{2}$ is not in the $\varepsilon -$%
galaxy of a rational standard, otherwise, $x_{n+1}$ and $\xi $ does not
become two successive elements of $S$. Hence, according to the lemma 2.5%
\begin{equation}
\left\{ 
\begin{array}{l}
\dfrac{x_{n+1}+\xi }{2}=\dfrac{s}{l}-\varepsilon \phi _{4} \\ 
\varepsilon l\cong 0\text{ },\text{ }l\cong +\infty \text{ },\text{ }\phi
_{4}\geq 0\text{ and }\phi _{4}\cong 0%
\end{array}%
\text{.}\right. \quad  \tag{2.12}
\end{equation}%
where $\dfrac{s}{l}$ is irreducible. Let $\overline{\gamma \text{ }}$ be an
unlimited natural number such that $\sqrt{\epsilon }.\overline{\gamma \text{ 
}}\cong 0$ and $\overline{N}=E\left( \dfrac{\overline{\gamma \text{ }}}{%
\sqrt{\epsilon }}\right) $. Let us take $\underline{N}=\max \left( \overline{%
N},k_{2},l\right) $. Then $\underline{N}\cong +\infty $ and is of the $%
E\left( \dfrac{\lambda }{\sqrt{\epsilon }}\right) $ form with $\lambda $ is
a positive unlimited real verifying $\sqrt{\epsilon }.\lambda \cong 0$. In
the other hand $\dfrac{h_{2}}{k_{2}}$ and $\dfrac{s}{l}$ are two elements of 
$\mathcal{F}_{\underline{N}}$ such that $\left[ \dfrac{h_{2}}{k_{2}}\text{ , 
}\dfrac{s}{l}\right] $ doesn't contain any rational standard and $k_{2}\cong
+\infty $ and $l\cong +\infty $. In this situation the lemma 2.4 is
applicable and consequently there is a finite sequence of irreducible
rational numbers $\left( \dfrac{s_{i}}{l_{i}}\right) _{1\leq i\leq e}$ such
that

\begin{equation*}
\dfrac{h_{2}}{k_{2}}=\dfrac{s_{1}}{l_{1}}<\dfrac{s_{2}}{l_{2}}<...<\dfrac{%
s_{e}}{l_{e}}=\dfrac{s}{l}
\end{equation*}%
where $e\cong +\infty $ and for $i=1,2,...,e-1$ we have :

\begin{equation*}
\dfrac{s_{i+1}}{l_{i+1}}-\dfrac{s_{i}}{l_{i}}=\varepsilon \phi \text{.}
\end{equation*}

\noindent Besides we have $\varepsilon l_{i}\cong 0$, $l_{i}\cong +\infty $
for $i=1,2,...,e$ ; $\dfrac{s_{e}}{l_{e}}-\dfrac{s_{1}}{l_{1}}=\varepsilon
\left( \dfrac{\omega }{2}+\phi _{4}-\phi _{2}\right) $.

In this paragraph we will associate to each $i\in \left\{ 1,2,...,e\right\} $
a vector $V_{i}$ in $\mathbb{Q}^{n+1}$ such that the $n$ first components of 
$V_{i}$ are in the $\varepsilon $-galaxie of the $n$ first components of $%
\left( x_{1},x_{2},...,x_{n},x_{n+1}\right) $, respectively. Whereas the $%
\left( n+1\right) -$th component of $V_{i}$ is equal to $\dfrac{s_{i}}{l_{i}}
$. Indeed, for $i=1$ apply lemma 2.2 to the system $\left(
l_{1}x_{1},l_{1}x_{2},...,l_{1}x_{n}\right) $ with the infinitesimal $%
\varepsilon l_{1}$:

\begin{equation*}
\left\{ 
\begin{array}{ccc}
\begin{array}{lllll}
l_{1}x_{i} & = & \dfrac{T_{i,1}}{t_{1}} & + & \left( \varepsilon
l_{1}\right) \pounds 
\end{array}
& ; & i=1,2,...,n \\ 
\epsilon l_{1}t_{1}\cong 0 &  & 
\end{array}%
\text{.}\right.
\end{equation*}

\noindent Hence $\left\{ 
\begin{array}{ccc}
\begin{array}{lllll}
x_{i} & = & \dfrac{T_{i,1}}{l_{1}t_{1}} & + & \varepsilon \pounds 
\end{array}
& ; & i=1,2,...,n \\ 
\varepsilon l_{1}t_{1}\cong 0 &  & 
\end{array}%
\text{.}\right. $ Then

\begin{equation}
\left\{ 
\begin{array}{ccc}
\begin{array}{lllll}
x_{i} & = & \dfrac{T_{i,1}}{l_{1}t_{1}} & + & \varepsilon \pounds 
\end{array}
& ; & i=1,2,...,n \\ 
\begin{array}{lllll}
x_{n+1} & = & \dfrac{T_{n+1,1}}{l_{1}t_{1}} & - & \lambda _{1}%
\end{array}
&  & 
\end{array}%
\right.  \tag{2.13}
\end{equation}%
where$T_{n+1,1}=s_{1}t_{1}$, $\varepsilon l_{1}t_{1}\cong 0$ and $\lambda
_{1}=\varepsilon \phi _{2}$. Then we obtain the vector $V_{1}=\left( \dfrac{%
T_{1,1}}{l_{1}t_{1}}\text{ },\text{ }\dfrac{T_{2,1}}{l_{1}t_{1}}\text{ },%
\text{ }...,\dfrac{T_{n+1,1}}{l_{1}t_{1}}\right) ^{T}$, where $x_{n+1}=%
\dfrac{T_{n+1,1}}{l_{1}t_{1}}=\dfrac{s_{1}}{l_{1}}$.

\noindent Again the application of the lemma 2.2 to the system $\left(
l_{2}x_{1},l_{2}x_{2},...,l_{2}x_{n}\right) $ with the infinitesimal $%
\varepsilon l_{2}$, gives:

$\left\{ 
\begin{array}{ccc}
\begin{array}{lllll}
l_{2}x_{i} & = & \dfrac{T_{i,2}}{t_{2}} & + & \left( \varepsilon
l_{2}\right) \pounds 
\end{array}
& ; & i=1,2,...,n \\ 
\epsilon l_{2}t_{2}\cong 0 &  & 
\end{array}%
\text{.}\right. $

\noindent Hence $\left\{ 
\begin{array}{ccc}
\begin{array}{lllll}
x_{i} & = & \dfrac{T_{i,2}}{l_{2}t_{2}} & + & \varepsilon \pounds 
\end{array}
& ; & i=1,2,...,n \\ 
\varepsilon l_{2}t_{2}\cong 0 &  & 
\end{array}%
\right. $. Then

\begin{equation}
\left\{ 
\begin{array}{ccc}
\begin{array}{lllll}
x_{i} & = & \dfrac{T_{i,2}}{l_{2}t_{2}} & + & \varepsilon \pounds 
\end{array}
& ; & i=1,2,...,n \\ 
\begin{array}{lllll}
x_{n+1} & = & \dfrac{T_{n+1,2}}{l_{2}t_{2}} & - & \lambda _{2}%
\end{array}
&  & 
\end{array}%
\right.  \tag{2.14}
\end{equation}

\noindent where $T_{n+1,2}=s_{2}t_{2}$ , $\varepsilon l_{2}t_{2}\cong 0$ and 
$\lambda _{2}=\varepsilon \phi _{2}+\varepsilon \phi $ with $0<\lambda
_{1}<\lambda _{2}$. Then we obtain the vector $V_{2}=\left( \dfrac{T_{1,2}}{%
l_{2}t_{2}}\text{ },\text{ }\dfrac{T_{2,2}}{l_{2}t_{2}}\text{ },\text{ }...,%
\dfrac{T_{n+1,2}}{l_{2}t_{2}}\right) ^{T}$, where $x_{n+1}=\dfrac{T_{n+1,2}}{%
l_{2}t_{2}}=\dfrac{s_{2}}{l_{2}}$.

\noindent Thus we construct the following vectors:

\begin{equation}
V_{i}=\left( \dfrac{T_{1,i}}{l_{i}t_{i}}\text{ },\text{ }\dfrac{T_{2,i}}{%
l_{i}t_{i}}\text{ },\text{ }...,\dfrac{T_{n+1,i}}{l_{i}t_{i}}\right)
^{T};\quad i=1,2,...,e  \tag{2.15}
\end{equation}%
where for $i=1,2,...,e$ : $x_{n+1}=\dfrac{T_{n+1,i}}{l_{i}t_{i}}-\lambda
_{i}=\dfrac{s_{i}}{l_{i}}-\lambda _{i}$ with $\varepsilon l_{i}t_{i}\cong 0$.

\noindent Besides $0<\varepsilon \phi _{2}=\lambda _{1}<\lambda
_{2}<...<\lambda _{e}=\dfrac{\varepsilon \omega }{2}+\varepsilon \phi _{4}$
and for $i=1,2,...,e-1$:

\begin{equation*}
\lambda _{i+1}-\lambda _{i}=\varepsilon \phi \text{.}
\end{equation*}

Let $h$ be the smallest integer such that $hNq\geq \underset{i}{\max }\left(
l_{i}t_{i}\right) $, then $\varepsilon hNq\cong 0$. On the other hand and
according to Robinson's lemma it exists an integer $W\cong +\infty $ such
that:%
\begin{equation*}
\varepsilon WhNq\cong 0\text{.}
\end{equation*}%
Put $K=hNq$. From $\left( 2.10\right) $:

\begin{equation}
\left\{ 
\begin{array}{ccc}
\begin{array}{lllllllll}
x_{i} & = & \dfrac{hNp_{i}}{hNq} & + & \varepsilon \phi & = & \dfrac{H_{i}}{K%
} & + & \varepsilon \phi%
\end{array}
& ; & i=1,2,...,n \\ 
\begin{array}{lllllllll}
x_{n+1} & = & \dfrac{hM}{hNq} & + & \varepsilon a & = & \dfrac{H_{n+1}}{K} & 
+ & \varepsilon a%
\end{array}
&  & 
\end{array}%
\right.  \tag{2.16}
\end{equation}

\noindent where $\varepsilon K\cong 0$, $K\geq \underset{i}{\max }\left(
l_{i}t_{i}\right) $.

\noindent Let $\overline{W}=\min \left( W\text{, }\dfrac{\omega }{2}+\phi
_{4}-\phi _{2}\right) $ and $\dfrac{T_{n+1,i_{0}}}{l_{i_{0}}t_{i_{0}}}$ be
the element of the sequence $\left( \dfrac{T_{n+1,i}}{l_{i}t_{i}}\right)
_{i=1,2,...,e}$ which is the farthest from $\dfrac{T_{n+1,1}}{l_{1}t_{1}}$
verifying

\begin{equation*}
\dfrac{T_{n+1,i_{0}}}{l_{i_{0}}t_{i_{0}}}-\dfrac{T_{n+1,1}}{l_{1}t_{1}}%
=\varepsilon \overline{\overline{W}}
\end{equation*}%
with $\overline{\overline{W}}\leq \overline{W}$. One notices that $\overline{%
\overline{W}}\cong +\infty $ because by construction $\overline{W}-\overline{%
\overline{W}}=\phi $.

\noindent Let $R\geq 1$ be the integer such that $Rl_{i_{0}}t_{i_{0}}\leq
K<\left( R+1\right) l_{i_{0}}t_{i_{0}}$. In this case $Rl_{i_{0}}t_{i_{0}}$
and $K$ are of the same order of magnitude i.e. : $\dfrac{K}{%
Rl_{i_{0}}t_{i_{0}}}=\delta $ where $\delta $ is a positive appreciable.
Consider, the rationals of the following vector:

\begin{equation}
\left( \dfrac{RT_{1,i_{0}}}{Rl_{i_{0}}t_{i_{0}}}\text{, }\dfrac{RT_{2,i_{0}}%
}{Rl_{i_{0}}t_{i_{0}}}\text{, ..., }\dfrac{RT_{n,i_{0}}}{Rl_{i_{0}}t_{i_{0}}}%
\text{, }\dfrac{RT_{n+1,i_{0}}}{Rl_{i_{0}}t_{i_{0}}}\right) ^{T}\text{.} 
\tag{2.17}
\end{equation}%
Where the $n$ first components of $\left(
x_{1},x_{2},...,x_{n},x_{n+1}\right) $ are in the $\varepsilon -$galaxies of
the $n$ first components of the (2.17), respectively. Whereas $x_{n+1}$ is
far from the last component of (2.17) by $\varepsilon \overline{\overline{W}}%
+\varepsilon \phi _{2}$. We will search a positive integer $j_{0}$ for which
the rational $\dfrac{RT_{n+1,i_{0}}+j_{0}H_{n+1}}{Rl_{i_{0}}t_{i_{0}}+j_{0}K}
$ becomes equal to $\dfrac{H_{n+1}}{K}+\varepsilon a+\varepsilon \phi $ i.e.
equal to $x_{n+1}+\varepsilon \phi $. Indeed, put

\begin{equation}
\Delta _{j}=\dfrac{RT_{n+1,i_{0}}+jH_{n+1}}{Rl_{i_{0}}t_{i_{0}}+jK}-\dfrac{%
H_{n+1}}{K}\text{.}  \tag{2.18}
\end{equation}%
Then $\Delta _{j}=\dfrac{\Delta }{1+j\delta }$ where $\Delta $ is the
distance between $\dfrac{RT_{n+1,i_{0}}}{Rl_{i_{0}}t_{i_{0}}}$ and $\dfrac{%
H_{n+1}}{K}$ which is equal to $\varepsilon \overline{\overline{W}}%
+\varepsilon \phi _{2}+\varepsilon a$.

\noindent Put $\dfrac{\Delta }{1+j\delta }=\varepsilon a$. For this $%
1+j\delta =\dfrac{\Delta }{\varepsilon a}$. Hence

\begin{equation*}
\begin{array}{lllll}
j & = & \dfrac{1}{\delta }\left( \dfrac{\Delta -\varepsilon a}{\varepsilon a}%
\right) & = & \dfrac{1}{\delta }\left( \dfrac{\varepsilon \overline{%
\overline{W}}+\varepsilon \phi _{2}}{\varepsilon a}\right) \\ 
&  &  & = & \dfrac{\overline{\overline{W}}+\phi _{2}}{\delta a}\cong +\infty%
\end{array}%
\text{. }
\end{equation*}%
Let us take $j_{0}=E\left( \dfrac{\overline{\overline{W}}+\phi _{2}}{\delta a%
}\right) $, hence $j_{0}=\dfrac{\overline{\overline{W}}+\phi _{2}}{\delta a}%
-\rho $ with $\rho \in \left[ 0,1\right[ $. Then $\Delta _{j_{0}}=\dfrac{%
\Delta }{1+j_{0}\delta }$. After the substitution by the value of $\Delta $
and of $j_{0}$:

\begin{equation*}
\begin{array}{lll}
\Delta _{j_{0}} & = & a.\dfrac{\varepsilon \overline{\overline{W}}%
+\varepsilon \phi _{2}+\varepsilon a}{a+\overline{\overline{W}}+\phi
_{2}-\rho a\delta } \\ 
& = & \varepsilon a\left( \dfrac{\overline{\overline{W}}+\phi _{2}+a}{%
\overline{\overline{W}}+\phi _{2}+a-\rho a\delta }\right)%
\end{array}%
\text{. }
\end{equation*}

\noindent Hence

\begin{equation*}
\begin{array}{lll}
\Delta _{j_{0}} & = & \varepsilon a.\dfrac{\left( \overline{\overline{W}}%
+\phi _{2}+a\right) }{\left( \overline{\overline{W}}+\phi _{2}+a\right)
\left( 1-\dfrac{\rho a\delta }{\overline{\overline{W}}+\phi _{2}+a}\right) }
\\ 
& = & \varepsilon a\dfrac{1}{1-\phi }\quad \text{.}%
\end{array}%
\text{ }
\end{equation*}%
Since $\dfrac{1}{1-\phi }=1+\phi $, then :

\begin{equation}
\Delta _{j_{0}}=\varepsilon a+\varepsilon \phi \text{.}  \tag{2.19}
\end{equation}%
On the other hand $j_{0}$ and $\overline{\overline{W}}$ are of the same
order of magnitude; indeed:

\begin{equation*}
\begin{array}{lll}
\dfrac{j_{0}}{\left( \overline{\overline{W}}\right) } & = & \dfrac{1}{\left( 
\overline{\overline{W}}\right) }\left( \dfrac{\overline{\overline{W}}+\phi
_{2}-\rho a\delta }{a\delta }\right) \\ 
& = & \dfrac{1+\phi }{a\delta }\text{.}%
\end{array}%
\end{equation*}%
Therefore $\dfrac{j_{0}}{\left( \overline{\overline{W}}\right) }=A$ with $A$
is appreciable, hence $j_{0}=A\overline{\overline{W}}$. Since $\overline{%
\overline{W}}\leq \overline{W}$ one has: $j_{0}=A\overline{\overline{W}}\leq
A\overline{W}\leq AW$.\newline

\noindent \textbf{Lemma} \textbf{2.9.} The denominator of \textit{\ }$\dfrac{%
RT_{i,i_{0}}+j_{0}H_{i}}{Rl_{i_{0}}t_{i_{0}}+j_{0}K}$ ( $i=1,2,...,n+1$)
verifies\textit{\ }$\varepsilon \left( Rl_{i_{0}}t_{i_{0}}+j_{0}K\right)
\cong 0$\textit{\ }and for $i=1,2,...,n,n+1$\ we have:

\begin{equation}
x_{i}=\dfrac{RT_{i,i_{0}}+j_{0}H_{i}}{Rl_{i_{0}}t_{i_{0}}+j_{0}K}%
+\varepsilon \phi  \tag{2.20}
\end{equation}

\noindent \textbf{Proof.}

\begin{equation*}
\begin{array}{lll}
Rl_{i_{0}}t_{i_{0}}+j_{0}K & = & \dfrac{K}{\delta }+j_{0}K \\ 
& = & K\left( \dfrac{1}{\delta }+j_{0}\right) \text{.}%
\end{array}%
\end{equation*}%
Hence $Rl_{i_{0}}t_{i_{0}}+j_{0}K\leq K\left( \dfrac{1}{\delta }+AW\right) $%
. From the fact that $\varepsilon WK\cong 0$; $A$ and $\delta $ are two
appreciable numbers, we have $\varepsilon \left(
Rl_{i_{0}}t_{i_{0}}+j_{0}K\right) \cong 0$. On the other hand for $i=n+1$ we
have from (2.16) $x_{n+1}=\dfrac{H_{n+1}}{K}+\varepsilon a$ and from (2.18)
and (2.19)

\begin{equation*}
\Delta _{j_{0}}=\dfrac{RT_{n+1,i_{0}}+j_{0}H_{n+1}}{%
Rl_{i_{0}}t_{i_{0}}+j_{0}K}-\dfrac{H_{n+1}}{K}=\varepsilon a+\varepsilon
\phi \text{.}
\end{equation*}%
Hence $\dfrac{RT_{n+1,i_{0}}+j_{0}H_{n+1}}{Rl_{i_{0}}t_{i_{0}}+j_{0}K}%
-\varepsilon \phi =\dfrac{H_{n+1}}{K}+\varepsilon a=x_{n+1}$, this means that

\begin{equation*}
x_{n+1}=\dfrac{RT_{n+1,i_{0}}+j_{0}H_{n+1}}{Rl_{i_{0}}t_{i_{0}}+j_{0}K}%
+\varepsilon \phi \text{.}
\end{equation*}%
For $i=1,2,...,n$ \ we know from $(2.17)$ that:

\begin{equation}
\left\vert \dfrac{RT_{i,i_{0}}}{Rl_{i_{0}}t_{i_{0}}}-x_{i}\right\vert
=\varepsilon \pounds \text{.}  \tag{2.21}
\end{equation}

\noindent Hence

\begin{quote}
$%
\begin{tabular}{lll}
$\left| \dfrac{RT_{i,i_0}}{Rl_{i_0}t_{i_0}}-\dfrac{H_i}K\right| $ & $=$ & $%
\left| \dfrac{RT_{i,i_0}}{Rl_{i_0}t_{i_0}}-x_i+x_i-\dfrac{H_i}K\right| $ \\ 
& $\leq $ & $\left| \dfrac{RT_{i,i_0}}{Rl_{i_0}t_{i_0}}-x_i\right| +\left|
x_i-\dfrac{H_i}K\right| =\varepsilon \pounds +\varepsilon \phi =\varepsilon 
\pounds \quad $.%
\end{tabular}
$
\end{quote}

\noindent\ Therefore 
\begin{equation}
\left\vert \dfrac{RT_{i,i_{0}}}{Rl_{i_{0}}t_{i_{0}}}-\dfrac{H_{i}}{K}%
\right\vert =\left\vert \dfrac{RT_{i,i_{0}}K-H_{i}Rl_{i_{0}}t_{i_{0}}}{%
KRl_{i_{0}}t_{i_{0}}}\right\vert =\epsilon \pounds \text{.}  \tag{2.22}
\end{equation}

\noindent Then we have: $\left\vert \dfrac{RT_{i,i_{0}}+j_{0}H_{i}}{%
Rl_{i_{0}}t_{i_{0}}+j_{0}K}-\dfrac{H_{i}}{K}\right\vert =\left\vert \dfrac{%
RT_{i,i_{0}}K-H_{i}Rl_{i_{0}}t_{i_{0}}}{KRl_{i_{0}}t_{i_{0}}\left( 1+j_{0}.%
\dfrac{K}{Rl_{i_{0}}t_{i_{0}}}\right) }\right\vert =\dfrac{\epsilon \pounds 
}{1+j_{0}\delta }$. Since $j_{0}\cong +\infty $, then

\begin{equation*}
\left\vert \dfrac{RT_{i,i_{0}}+j_{0}H_{i}}{Rl_{i_{0}}t_{i_{0}}+j_{0}K}-%
\dfrac{H_{i}}{K}\right\vert =\varepsilon \phi \text{ }
\end{equation*}%
and seen that for $i=1,2,...,n$, the rational numbers $\dfrac{H_{i}}{K}$
are, respectively, in the $\varepsilon -$halos of $x_{1},x_{2},...,x_{n}$
then:

\begin{equation*}
x_{i}=\dfrac{RT_{i,i_{0}}+j_{0}H_{i}}{Rl_{i_{0}}t_{i_{0}}+j_{0}K}%
+\varepsilon \phi \text{.}
\end{equation*}%
So the lemma is proved.\newline

Since $\varepsilon \left( Rl_{i_{0}}t_{i_{0}}+j_{0}K\right) \cong 0$, then
if for $i=1,2,...,n,n+1$ one takes $\dfrac{RT_{i,i_{0}}+j_{0}H_{i}}{%
Rl_{i_{0}}t_{i_{0}}+j_{0}K}$ for $\dfrac{P_{i}}{Q}$ then

\begin{equation}
\left\{ 
\begin{array}{ccc}
x_{i} & = & \dfrac{P_{i}}{Q}+\varepsilon \phi \text{, }i=1,2,...,n+1 \\ 
& \varepsilon Q\cong 0 & 
\end{array}%
\text{.}\right.  \tag{2.23}
\end{equation}

\noindent In the case where $a<0$ we take $\xi $ the element of $S$ that
precedes $x_{n+1}$ i.e. $\xi <x_{n+1}$ ($S$ is ordered) and by doing, to a
symmetry near, as we did for the case $a>0$.

\noindent From $(2.8)$, $(2.9)$ and $(2.23)$ we have $A(n+1)$. Hence,
according to the external recurrence principle, the lemma 2.8 is proved. $%
\square $\newline

Let us return to the proof of theorem 2.1

\noindent Define for $Z=\left\{ x_{1},x_{2},...,x_{s}\right\} \subset \left[
0,1\right] $, the formula:

\begin{equation}
B\left( Z\right) ="\exists \left( \dfrac{P_{i}}{Q}\right) _{i=1,2,...,s}%
\text{such that :}\forall ^{st}m\in \mathbb{N}^{\ast }\text{ }G\left( Z\text{%
, }\left( \dfrac{P_{i}}{Q}\right) _{i=1,2,...,s}\text{, }m\right) " 
\tag{2.24}
\end{equation}%
where $G\left( Z\text{, }\left( \dfrac{P_{i}}{Q}\right) _{i=1,2,...,s}\text{%
, }m\right) \equiv \left\{ 
\begin{array}{cc}
\dfrac{1}{\varepsilon }\left\vert x_{i}-\dfrac{P_{i}}{Q}\right\vert \leq 
\dfrac{1}{m} & ;=1,2,...,s \\ 
\left\vert \varepsilon Q\right\vert \leq \dfrac{1}{m} & 
\end{array}%
\right. $ is internal.

\noindent Consider the set

\begin{equation}
L=\left\{ n\in \mathbb{N}^{\ast }:n\leq \left\vert S\right\vert \text{ \& }%
\forall s\in \left\{ 1\text{,}...\text{,}n\right\} \text{ }\forall Z=\left\{
x_{1},x_{2},...,x_{s}\right\} \subset S\text{ }:B(Z)\right\} \text{.} 
\tag{2.25}
\end{equation}%
where $S$\ is the set that has been constructed in the lemma 2.6 .Then

\begin{equation*}
L=\left\{ 
\begin{array}{c}
n\in \mathbb{N}^{\ast }:n\leq \left\vert S\right\vert \text{ \&}\forall s\in
\left\{ 1\text{,}...\text{,}n\right\} \text{ }\forall Z=\left\{
x_{1},x_{2},...,x_{s}\right\} \subset S\text{,} \\ 
\exists \left( \dfrac{P_{i}}{Q}\right) _{i=1,2,...,s}\text{ }\forall
^{st}m\in \mathbb{N}^{\ast }G\left( Z\text{, }\left( \dfrac{P_{i}}{Q}\right)
_{i=1,2,...,s}\text{, }m\right)%
\end{array}%
\right\} \text{.}
\end{equation*}

According to lemma 2.8, $\ L\supset $ $\left( \mathbb{N}^{\ast }\right)
^{\sigma }$. If $L$ is internal then, according to the Cauchy principle, it
must contain $\left( \mathbb{N}^{\ast }\right) ^{\sigma }$ strictly and
therefore there is an integer $\omega \cong +\infty $ and $\omega \in L$. If 
$L$ is external then by the idealization principle (I) we can write $L$ as
follows:\newline

\begin{equation*}
L=\left\{ 
\begin{array}{c}
n\in \mathbb{N}^{\ast }:n\leq \left\vert S\right\vert \text{ \&}\forall s\in
\left\{ 1\text{,}...\text{,}n\right\} \text{ }\forall Z=\left\{
x_{1},x_{2},...,x_{s}\right\} \subset S\text{,} \\ 
\forall ^{stfini}\text{ }M\text{ }\exists \left( \dfrac{P_{i}}{Q}\right)
_{i=1,2,...,s}\forall m\in M\text{ }G\left( Z\text{, }\left( \dfrac{P_{i}}{Q}%
\right) _{i=1,2,...,s}\text{, }m\right)%
\end{array}%
\right\} \text{.}
\end{equation*}%
where $M$ belongs to the set of finite parts of $\mathbb{N}^{\ast }$.
Therefore, $L$ is an halo ($\left[ 4\right] $, $\left[ 6\right] $). Of the
fact that $\left( \mathbb{N}^{\ast }\right) ^{\sigma }\subset L$ and no halo
is a galaxy (Fehrele principle), then $\left( \mathbb{N}^{\ast }\right)
^{\sigma }\underset{\neq }{\subset }L$. Hence it exists an integer $\omega
\cong +\infty $ and $\omega \in L$.

Consequently in the two cases ($L$ internal or external ) we finds that it
exists an integer $\omega \cong +\infty $ and $\omega \in L$, this signifies
that $\omega \leq \left\vert S\right\vert $.

By lemma 2.3, there is a finite part $F\subset \left[ 0,1\right] $
containing all standard elements of $\left[ 0,1\right] $ such that $%
\left\vert F\right\vert =\omega ^{\prime }\cong +\infty $ and $\omega
^{\prime }<\omega $. Then $F\cap S$ is a finite part of $S$ containing all
standard elements of $\left[ 0,1\right] $ with $\left\vert F\cap
S\right\vert \leq \left\vert F\right\vert =\omega ^{\prime }<\omega $. Put $%
F\cap S=\left\{ x_{1},x_{2},...,x_{n_{0}}\right\} $. Then $\exists \left( 
\dfrac{P_{i}}{Q}\right) _{i=1,2,...,n_{0}}$such that :$\left\{ 
\begin{array}{c}
x_{i}-\dfrac{P_{i}}{Q}=\varepsilon \phi \\ 
\varepsilon Q\cong 0\text{ };\text{ }i=1,2,...,n_{0}%
\end{array}%
\right. $. It follows that if $x\in \mathbb{R}$ is a standard then $%
\begin{tabular}{lll}
$x-E\left( x\right) $ & $=$ & $\dfrac{P_{i_{1}}}{Q}+\varepsilon \phi $%
\end{tabular}%
$ where $i_{1}\in \left\{ 1,2,...,n_{0}\right\} $ since $x-E\left( x\right) $
is a standard of $\left[ 0,1\right] $. Hence%
\begin{equation*}
\left\{ 
\begin{tabular}{lllll}
$x$ & $=$ & $E\left( x\right) +\dfrac{P_{i_{1}}}{Q}+\varepsilon \phi $ & $=$
& $\dfrac{P_{x}}{Q}+\varepsilon \phi $%
\end{tabular}%
\right.
\end{equation*}%
where $\varepsilon Q\cong 0$. Thus the proof is complete.\newline

\section{Deduction of the classical equivalent of the main result}

The theorem 2.1. can be written as follows

\begin{equation*}
\forall \varepsilon \left\{ \left( \forall ^{st}r\text{ }\left(
0<\varepsilon \leq r\right) \right) \Longrightarrow \exists q\text{ }\forall
^{st}x\text{ }\forall ^{st}t\text{ }\left( \parallel qx\parallel
<\varepsilon qt\text{ }\&\text{ }\varepsilon q\leq t\right) \right\}
\end{equation*}%
where $\varepsilon $, $r\in \mathbb{R}^{\ast +}$, $q\in \mathbb{N}$, $x\in 
\mathbb{R}$ and $t\in \mathbb{R}^{\ast +}$. By using the idealization
principle (I), the last formula is equivalent to

\begin{equation*}
\forall \varepsilon \left\{ \left( \forall ^{st}r\text{ }\left(
0<\varepsilon \leq r\right) \right) \Longrightarrow \forall ^{st\text{ }%
fini}X\text{ }\exists q\text{ }\forall \left( x\text{, }t\right) \in X\text{ 
}\left( \parallel qx\parallel <\varepsilon qt\text{ }\&\text{ }\varepsilon
q\leq t\right) \right\}
\end{equation*}%
where $X$ belongs to the set of finite parts of $\mathbb{R}\times \mathbb{R}%
^{\ast +}$. This last formula is equivalent to

\begin{equation*}
\forall ^{st\text{ }fini}X\forall \varepsilon \exists ^{st}r\left\{ \left(
0<\varepsilon \leq r\right) \Longrightarrow \exists q\text{ }\forall \left( x%
\text{, }t\right) \in X\text{ }\left( \parallel qx\parallel <\varepsilon qt%
\text{ }\&\text{ }\varepsilon q\leq t\right) \right\} \text{.}
\end{equation*}%
Again, by using the idealization principle (I), the last formula is
equivalent to

\begin{equation*}
\forall ^{st\text{ }fini}X\text{ }\exists ^{st\text{ }fini\text{ }}R\text{ }%
\forall \varepsilon \text{ }\exists r\in R\text{ }\left\{ \left(
0<\varepsilon \leq r\right) \Longrightarrow \exists q\text{ }\forall \left( x%
\text{, }t\right) \in X\text{ }\left( \parallel qx\parallel <\varepsilon qt%
\text{ }\&\text{ }\varepsilon q\leq t\right) \right\} \text{.}
\end{equation*}%
where $R$ belongs to the set of finite parts of $\mathbb{R}^{\ast +}$. By
the transfer principle (T), this last formula is equivalent to

\begin{equation*}
\forall ^{fini}X\text{ }\exists ^{fini\text{ }}R\text{ }\forall \varepsilon 
\text{ }\exists r\in R\text{ }\left\{ \left( 0<\varepsilon \leq r\right)
\Longrightarrow \exists q\text{ }\forall \left( x\text{, }t\right) \in X%
\text{ }\left( \parallel qx\parallel <\varepsilon qt\text{ }\&\text{ }%
\varepsilon q\leq t\right) \right\} \text{.}
\end{equation*}

This last formula is exactly the main theorem announced in the abstract.
Indeed, if$\ X=\left\{ \left( x_{1}\text{, }t_{1}\right) \text{, }\left(
x_{2}\text{, }t_{2}\right) \text{, ..., }\left( x_{n}\text{, }t_{n}\right)
\right\} $ is a finite part of $\mathbb{R}\times \mathbb{R}^{\ast +}$, then
there exist a finite part $R$ of $\mathbb{R}^{\ast +}$ such that for all $%
\varepsilon >0$ there exists $r\in R$ such that if $0<\varepsilon \leq r$
then there exist rational numbers $\left( \dfrac{p_{i}}{q}\right)
_{i=1,2,...,n}$ such that:

\begin{equation*}
\left\{ 
\begin{array}{c}
\left\vert x_{i}-\dfrac{p_{i}}{q}\right\vert \leq \varepsilon t \\ 
\varepsilon q\leq t%
\end{array}%
\right\vert ,=1,2,...,n\text{.}
\end{equation*}

\begin{equation*}
References
\end{equation*}

\noindent $\left[ 1\right] $ A. BOUDAOUD, \textit{Mod\'{e}lisation de ph\'{e}%
nom\`{e}nes discrets et approximations diophantiennes infinit\'{e}simales},
Maghreb Mathematical Review Vol 1$\left( 1\right) $, June 1992.

\noindent $\left[ 2\right] $ F.Diener, G.Reeb, \textit{Analyse Non Standard}%
, Hermann \'{e}diteurs des sciences et des arts . 1989 .

\noindent $\left[ 3\right] $ G.H. HARDY \& E.M. WRIGHT, \textit{An
introduction to the theory of numbers}, Clarendon Press, Oxford, 1979.

\noindent $\left[ 4\right] $ F. Koudjeti, \textit{Elements of External
Calculus with an application to mathematical Finance, }PhD thesis, Groninjen
(The Netherlands), 1995.

\noindent $\left[ 5\right] $ R. Lutz, M.Goze, \textit{Non Standard Analysis.
A practical guide with applications}, Lecture note in Math. N$^{0}$881,
Springer Verlag (1981).

\noindent $\left[ 6\right] $ E. Nelson, Internal set theory : \textit{A new
approach to non standard analysis}, bull. Amer. Math. Soc. 83(1977)
1165-1198.

\noindent $\left[ 7\right] $ W.M. SCHMIDT, \textit{Diophantine approximation}%
, Lectures Notes in Mathematics, N$^{0}$785, Springer-Verlag, berlin, (1980).

\noindent $\left[ 8\right] $ I. P. Van den Berg, \textit{Nonstandard
Asymptotic Analysis, }Volume 1249 of Lecture Notes in Mathematics. Springer
Verlag, 1987.

\newpage

\end{document}